\documentclass{amsart}
\usepackage{color}

\usepackage{amsmath} 
\usepackage{amssymb}
\usepackage{mathrsfs}

\newtheorem{theorem}{Theorem}[section] 
\newtheorem{claim}[theorem]{Claim}

\newtheorem{conclusion}[theorem]{Conclusion}
\newtheorem{observation}[theorem]{Observation}

\theoremstyle{definition}
\newtheorem{definition}[theorem]{Definition}

\newtheorem{fact}[theorem]{Fact}
\newtheorem{discussion}[theorem]{Discussion}

\theoremstyle{remark}
\newtheorem{remark}[theorem]{Remark}
\newtheorem{question}[theorem]{Question}
\newtheorem{notation}[theorem]{Notation}

\newcommand{\inv}{{\rm inv}}

\newcommand{\tf}{{\rm tf}}
\newcommand{\trf}{{\rm trf}}

\newcommand{\Mod}{{\rm Mod}}

\newcommand{\Kernel}{{\rm Kernel}}
\newcommand{\otp}{{\rm otp}}

\newcommand{\ctrf}{{\rm ctrf}}
\newcommand{\crtf}{{\rm crtf}}
\newcommand{\univ}{{\rm univ}}

\newcommand{\pp}{{\rm pp}}

\newcommand{\bd}{{\rm bd}}

\newcommand{\nacc}{{\rm nacc}}

\newcommand{\rs}{{\rm rs}}

\newcommand{\LST}{{\rm LST}}

\newcommand{\cov}{{\rm cov}}

\newcommand{\rtf}{{\rm rtf}}

\newcommand{\Min}{{\rm Min}}

\newcommand{\rest}{{\restriction}}

\newcommand{\wilog}{{\rm without loss of generality}}
\newcommand{\Wilog}{{\rm Without loss of generality}}

\newcommand{\then}{{\underline{then}}}
\newcommand{\when}{{\underline{when}}}
\newcommand{\Then}{{\underline{Then}}}

\newcommand{\If}{{\underline{if}}}
\newcommand{\Iff}{{\underline{iff}}}
\newcommand{\mn}{{\medskip\noindent}}
\newcommand{\sn}{{\smallskip\noindent}}

\newcommand{\cA}{{\mathscr A}}

\newcommand{\cB}{{\mathscr B}}

\newcommand{\bbZ}{{\mathbb Z}}

\newcommand{\cH}{{\mathscr H}}

\newcommand{\bbL}{{\mathbb L}}
\newcommand{\bbN}{{\mathbb N}}
\newcommand{\cM}{{\mathscr M}}

\newcommand{\gk}{{\mathfrak k}}
\newcommand{\gB}{{\mathfrak B}}
\newcommand{\gC}{{\mathfrak C}}
\newcommand{\gK}{{\mathfrak K}}

\newcommand{\bbP}{{\mathbb P}}
\newcommand{\cP}{{\mathscr P}}

\newcommand{\gT}{{\mathfrak T}}

\newcommand{\bbQ}{{\mathbb Q}}

\newcommand{\cT}{{\mathscr T}}

\newcommand{\cf}{{\rm cf}}
\newcommand{\pr}{{\rm pr}}

\newcount\skewfactor
\def\mathunderaccent#1#2 {\let\theaccent#1\skewfactor#2
\mathpalette\putaccentunder}
\def\putaccentunder#1#2{\oalign{$#1#2$\crcr\hidewidth
\vbox to.2ex{\hbox{$#1\skew\skewfactor\theaccent{}$}\vss}\hidewidth}}

\newenvironment{PROOF}[2][\proofname.]
   {\begin{proof}[#1]\renewcommand{\qedsymbol}{\textsquare$_{\rm #2}$}}
   {\end{proof}}

\usepackage[hidelinks]{hyperref}

\begin{document}
\makeatletter\def\shfiuwefootnote{\gdef\@thefnmark{}\@footnotetext}\makeatother\shfiuwefootnote{Version 2017-04-06\_12. See \url{https://shelah.logic.at/papers/820/} for possible updates.}

\title {Universal Structures \\
 Sh820}
\author {Saharon Shelah}
\address{Einstein Institute of Mathematics\\
Edmond J. Safra Campus, Givat Ram\\
The Hebrew University of Jerusalem\\
Jerusalem, 91904, Israel\\
 and \\
 Department of Mathematics\\
 Hill Center - Busch Campus \\ 
 Rutgers, The State University of New Jersey \\
 110 Frelinghuysen Road \\
 Piscataway, NJ 08854-8019 USA}
\email{shelah@math.huji.ac.il}
\urladdr{http://shelah.logic.at}
\thanks{This research was partially
supported by the German-Israel Foundation for Scientific Research and
Development.
I would like to thank Alice Leonhardt for the beautiful typing.
First Typed - 98/June/10}



 

\date{April 6, 2017}

\begin{abstract}
We deal with the existence of universal members in a given
 cardinality for several classes.  First we deal with classes  
of Abelian groups,  specifically with the existence of universal 
members in cardinalities
 which are strong limit singular of countable cofinality or $\lambda =
 \lambda^{\aleph_0}$.  We use
 versions of being reduced replacing $\bbQ$ by a subring (defined by a sequence
 $\bar t$) and get quite accurate results for existence of universal
 in a cardinal, for embedding and for pure embeddings.
Second, we deal with (variants of) the oak property (from a work 
of D\v{z}amonja and the author),
 a property of complete first order theories, sufficient for the
 non-existence of universal models under suitable cardinal
 assumptions.  Third, we prove that the oak property 
holds for the class of groups (naturally interpreted, so for
 quantifier free formulas) and deal more with the existence of universals.
\end{abstract}

\maketitle
\numberwithin{equation}{section}
\setcounter{section}{-1}
\newpage

\centerline {Annotated Content}
\bigskip

\noindent
\S0 \quad Introduction, \pageref{intro}
\bigskip

\noindent
\S1 \quad More on Abelian groups, \pageref{moreon}
\mn
\begin{enumerate}
\item[${{}}$]   [We say more on some classes of Abelian groups; mainly
torsion free such that no non-zero $x$ is divisible by
$\prod\limits_{\ell < n} t_\ell$ for every $n$ (where $t_n \ge 2$).  We
get results on existence and non-existence of universal structures in 
cardinals like $\lambda = \lambda^{\aleph_0}$ and $\beth_\omega$, that
is, $\lambda = \sum\limits_{n} \lambda_n,\lambda_n =
(\lambda_n)^{\aleph_0}$.  For $\lambda = \lambda^{\aleph_0}$ we get
characterizations by $\bar t$, which is different for embedding and
pure embedding.  For strong limit $\lambda$ of cofinality $\aleph_0$, 
we use a general criterion for existence.]
\end{enumerate}
\bigskip

\noindent
\S2 \quad The class of groups, \pageref{theclass}
\mn
\begin{enumerate}
\item[${{}}$]   [We prove that the class of groups
has the oak property (from \cite{Sh:710}).]
\end{enumerate}
\bigskip

\noindent
\S3 \quad On the oak property, \pageref{ontheoak}
\mn
\begin{enumerate}
\item[${{}}$]   [We continue \cite{Sh:710}, deal with singular
  cardinals and a weaker relative of the property.]
\end{enumerate}
\newpage

\section {Introduction} \label{intro}

On the existence of universal structures see Kojman-Shelah \cite{Sh:409} and
history there, and a more recent survey D\v{z}amonja \cite{Dj05}.  Of
course, a complete first order theory $T$ has a universal model in $\lambda$
for ``elementary embeddings"  when
$\lambda = 2^{< \lambda} > |T|$; this is true also for similar classes,
i.e. for a.e.c. with amalgamation and the JEP and LST number $< \lambda$.  The
question which interests us is whether there are additional cases
(mainly for elementary classes and more generally a.e.c. as above).  
But here we deal with some
specific classes and embeddability notion.

Now \S1 deals mainly with Abelian groups; it continues Kojman-Shelah
\cite{Sh:455} and \cite{Sh:456}, \cite{Sh:552} and \cite{Sh:622}.  The 
second section deals with the class of groups; it continues 
Usvyatsov-Shelah \cite{Sh:789} but does not rely on it.  The third section
deals with the oak property continuing D\v{z}amonja-Shelah
\cite{Sh:710}, dealing with the case of singular cardinals.

The second section deals with the class of all groups, certainly an
important one.  Is this class complicated?  Under several yard-sticks
it certainly is: its first order theory is undecidable, etc., and it
has the quantifier-free order property (even the class of (universal) locally
finite groups, has this property, see Macintyre-Shelah 
\cite{Sh:55}) and by \cite{Sh:789} it
has the SOP$_3$ (3-strong order property).  But this does not
exclude positive answers for other interpretations.  By \cite{Sh:789}
it has the NSOP$_4$ (4-strong non-order property), however we do not know
much about this family of classes (though we have hopes).

A recent relevant work is \cite{Sh:1029}, \cite{Sh:F1425}, 
giving new sufficient
conditions for ``no universal", in particular for groups and hopefully
\cite{Sh:F1436} on classes of Abelian groups.

Here we consider the oak property from
D\v{z}amonja-Shelah \cite{Sh:710}, a relative of the tree property, (hence
the name).  We prove that 
the class of groups has the oak property, hence it follows that in some
cardinals it has no universal member.  

There is reasonable evidence for the class of linear orders
being complicated, practically maximal for the universal spectrum
problem, see \cite{Sh:409}.  

So a specific conclusion is:

\begin{conclusion}
\label{y3}
1) The class of groups has the oak property, see Definition \ref{b3}.

\noindent
2) If $\lambda$ satisfies, e.g., $(*)$ below \then \, there is no
universal group of cardinality $\lambda$ when:
\mn
\begin{enumerate}
\item[$(*)$]   $(a) \quad \kappa = \cf(\mu) < \mu$
\sn
\item[${{}}$]    $(b) \quad \lambda = \mu^{++} < \pp_{J^{\bd}_\kappa}(\mu)$
\sn
\item[${{}}$]   $(c) \quad \alpha < \mu \Rightarrow |\alpha|^\kappa < \mu$.
\end{enumerate}
\end{conclusion}

\begin{PROOF}{\ref{y3}}
1) By \ref{b6}.

\noindent
2) By part (1) and \cite{Sh:710}, more exactly by \ref{c3}.
\end{PROOF}

In \S3 we deal with the oak property per se, continuing \cite{Sh:710},
showing non-existence of universal in singular cardinals and dealing
with a weaker relative, the weak oak property.

Concerning the first section note that strong limit singular cardinal
$\lambda$ is a case where it is easier to have a universal model,
particularly when $\lambda$ has cofinality $\aleph_0$.

\noindent
So the canonical case seems to be $\beth_\omega$.  Examples of such 
positive (= existence) results are
\mn
\begin{enumerate}
\item[$(a)$]  \cite[Th.3.1,p.266]{Sh:26}, where it is proved that: 

if $\lambda$ is strong limit singular, \then  

$\{G:G \text{ a graph with } \le \lambda \text{ nodes each of valency }
< \lambda\}$ has a universal 

member under embedding onto induced subgraphs
\sn
\item[$(b)$]    Grossberg-Shelah \cite[\S5,Corollary 27,pgs.301-302]{Sh:174}:
\sn
\begin{enumerate}
\item[$(\alpha)$]   if $\lambda$ is ``large enough" then
similar results hold for quite general classes (e.g. locally finite groups) 
\underline{where} large enough means: $\lambda$ (is strong limit, of 
cofinality $\aleph_0$ and) is above a compact cardinal (which is quite large).

More specifically, 
\sn
\item[$(\beta)$]  if $\mu$ is strong limit of cofinality $\aleph_0$
above a compact cardinal $\kappa$ and, e.g., the class ${\gK}$ is
the class of models of $T \subseteq \bbL_{\kappa,\aleph_0},|T| <
\mu$ partially ordered by $\prec_{\bbL_{\kappa,\omega}}$, \then \, we
can split ${\gK}$ to $\le 2^{|T|+\kappa}$ classes each having a
universal model of cardinality $\mu$ under
$\prec_{\bbL_{\kappa,\aleph_0}}$-embeddings.   
\end{enumerate}
\end{enumerate}
\mn
See more in \cite{Sh:1011} generalizing the so-called special models.
Claim \ref{a43} below continues this, i.e., it deals with strong limit
cardinal $\mu > \cf(\mu) = \aleph_0$, compared to \cite{Sh:174} 
omitting the set
theoretic assumption on compact cardinal at the expense of
strengthening the model theoretic assumption.

There are natural examples where this can be applied; e.g. the class
of torsion free Abelian groups $G$ which are reduced (i.e., we
cannot embed the rational into $G$), \underline{but} the order is $G_1
\le_{\langle n!:n<\omega\rangle} G_2$ which means $G_1 \subseteq G_2$ but
$G_1$ is closed inside $G_2$ under the $\bbZ$-adic metric; so
also $G_2/G_1$ is reduced.  The application of \ref{a43} 
to such classes is in \ref{a37}(1)(2).
Earlier in \ref{a6} we prove related positive results for the
easier cases of complete members (for $\lambda$ satisfying $\lambda =
\lambda^{\aleph_0}$ or $\lambda$ the limit of such cardinals).

We also get some negative results, i.e., non-existence of universal
members in \ref{a18}(2), \ref{a28}.  We deal more generally with
$K^{\rtf}_{\bar t}$, the reduced torsion free Abelian group $G$ such
that for no $x \in G,x \ne 0$ and $x$ is divisible by $t_{<n} =
\prod\limits_{\ell < n} t_\ell$ for every $n$.  We sort out the existence of
universal members of cardinality $\lambda = \lambda^{\aleph_0}$ for
$K^{\rtf}_{\bar t,\lambda}$ under embeddings and under pure
embeddings, getting complete (but different) answers for $\lambda =
\lambda^{\aleph_0}$. 

Recall that classes of Abelian groups are related to the classes of trees
with $\omega +1$ levels.  The parallel of ``Abelian groups under pure
embedding" is the case of such trees, in fact, non-existence of
universals for Abelian groups under pure embedding implies the 
non-existence of such universal trees.

We thank the referee for many helpful comments.

\begin{notation}
\label{z2}a
1) For a set $A,|A|$ is its cardinality but for a structure $M$ its
   cardinality is $\|M\|$ while its universe is $|M|$; this apply
   e.g. to groups.

\noindent
2) $\bar t$ will denote an $\omega$-sequence of natural numbers $\ge 2$.

\noindent
3) We use $G,H$ for groups, $M,N$ for general models.

\noindent
4) Let $\gk$ denote a pair $(K_{\gk},\le_{\gk})$, we may 
say a class $\gk$ instead of a pair, where:
\mn
\begin{enumerate}
\item[$(a)$]  $K_{\gk}$ is a class of $\tau_{\gk}$-structures
\sn
\item[$(b)$]  $\le_{\gk}$ is a partial order on $K_{\gk}$ such that $M
  \le_{\gk} N \Rightarrow M \subseteq N$
\sn
\item[$(c)$]  both $K_{\gk}$ and $\le_{\gk}$ are closed under
  isomorphisms.
\end{enumerate}
\mn
4A) We say $f:M \rightarrow N$ is a $\le_{\gk}$-embedding when $f$ is
an isomorphism from $M$ onto some $M_1 \le_{\gk} N$.

\noindent
5) If $T$ is a first order theory then $\Mod_T$ is the pair
$(\text{mod}_T,\le_T)$ where $\text{ mod}_T$ is the class of models of $T$ and
   $\le_T$ is: $\prec$ if $T$ is complete, $\subseteq$ if $T$ is not
   complete.

\noindent
6) We may write $T$ instead of $\Mod_T$, e.g. in Definition \ref{z3} below.
\end{notation}

\begin{definition}
\label{z3}
1) For a class $\gk$ and a cardinal $\lambda$, a set $\{M_i:i<i^*\}$
of models from $\gk$, are \underline{jointly universal} for $\lambda$  
\when \, for every $N \in K_{\gk}$ of size $\lambda$,
there is an $i<i^*$ and an $\le_{\gk}$-embedding of $N$ into $M_i$.

\noindent
2) For $\gk$ and $\lambda$ as above, let (if $\mu = \lambda$ we may
   omit $\mu$)

\begin{equation*}
\begin{array}{clcr}
\univ(\gk,\mu,\lambda) := &\min\{|\cM|:\cM \text{ is a family of
  members of } K_{\gk} \text{ each} \\
  &\text{ of cardinality } \le \mu \text{ which is jointly universal} \\
  &\text{ for models of } \gk \text{ of size } \lambda\}.
\end{array}
\end{equation*}
\end{definition}

\begin{remark}
To help understanding Definition \ref{z3}, note that $\univ(T,\lambda)
= 1$ iff there is a universal model of $T$ of size $\lambda$.  Note
that some of the classes we consider are not abstract elementary
classes.  Some have ``weak failure" say $\bbZ$-adically complete
free Abelian free groups which are torsion free, if $M_n \le M_{n+1}$ then
$\bigcup\limits_{n} M_n$ is not necessarily complete.  We can take a
completion; more seriously for some $\gk$ and the 
$M_n$'s there are contradictory completions. 
\end{remark}

\noindent
Recall
\begin{definition}
\label{z18}
For an ideal $J$ on a set $A$ and a set $B$ let $\mathbf U_J(B) = 
\Min\{|{\cP}|:{\cP}$ is a family of subsets of $B$, each of 
cardinality $\le |A|$ such that 
for every function $f$ from $A$ into $B$ for some 
$u \in {\cP}$ we have $\{a \in A:f(a) \in u\} \in J^+\}$.  Clearly 
only $|B|$ matters so we normally write
$\mathbf U_J(\lambda)$ (see on it \cite{Sh:589}).
\end{definition}
\newpage

\section {More on Abelian groups} \label{moreon}

Earlier versions of this 
section originally was part of \cite{Sh:622} and earlier of
\cite{Sh:552}, but as the papers were too long, it was delayed. 

\begin{remark}
\label{a3}
Inspite of all cases dealt with in \cite{Sh:552}, there 
are still some ``missing" cardinals (see
discussion in \cite[\S0]{Sh:622}).  Concerning $\lambda$ singular satisfying
$2^{\aleph_0} < \mu^+ < \lambda < \mu^{\aleph_0}$, clearly
\cite[2.8=2.7t,3.14=3.12t]{Sh:622}, \cite{Sh:g} indicates that at least 
for most such cardinals there is no universal: as if 
$\chi \in (\mu^+,\lambda)$ is regular, then 
$\cov(\lambda,\chi^+,\chi^+,\chi) < \mu^{\aleph_0}$.
\end{remark}

\noindent
Let us mention concerning positive results on 
Case 1 (from \cite[\S0]{Sh:622}), see Definition \ref{a9} below.
(See Fuchs \cite{Fu} on Abelian groups). 

\begin{claim}
\label{a6}
1) If $\lambda = \lambda^{\aleph_0}$ \then \,
in the class $(K^{\rtf}_\lambda,\le_{\pr})$, defined in \ref{a9}(5) below 
there is a universal member, in fact it is homogeneous universal. 

\noindent
2) If $\lambda = \sum\limits_{n < \omega} \lambda_n$ and $\aleph_0 \le
\lambda_n = (\lambda_n)^{\aleph_0} < \lambda_{n+1}$ 
\then \, in $({\gK}^{\rtf}_\lambda,\le_{\pr})$ there is a 
universal member (the parallel of special models for first order theories). 

\noindent
3) $(K^{\rtf},\le_{\pr})$ has the amalgamation and JEP; is an
a.e.c. (see \cite{Sh:h}) and is stable in $\lambda$ if $\lambda =
\lambda^{\aleph_0}$. 
\end{claim}

\noindent
We shall prove \ref{a6} below, but first
\begin{definition}
\label{a9}
1) $K^{\tf}_\lambda$ is the class of torsion-free Abelian 
groups of cardinality $\lambda$.  Let
$K^{\tf} = \cup\{K^{\tf}_\lambda:\lambda$ a cardinal$\}$ and similarly
$K^{\tf}_{\le \lambda}$.

\noindent
1A)  $K^{\rtf}_{\bar t,\lambda}$ is the class of $G \in 
K^{\tf}_\lambda$ such that there is no $x \in G \backslash \{0\}$
divisible by $\prod\limits_{\ell < k} t_\ell$ for every $k < \omega$
recalling \ref{z2}(2).

\noindent
1B) Let $K^{\rtf}_{\bar t} =
\cup\{K^{\rtf}_{\bar t,\lambda}:\lambda$ a cardinal$\}$.
 
\noindent
1C) $G \in K^{\rtf}_{\bar t}$ is called $\bar t$-complete \when \, every
Cauchy sequence under $d_{\bar t}$ in $G$ has a limit where $d_{\bar
  t}$ is defined in \ref{a9}(3) below.

\noindent
2)  Let
\mn
\begin{enumerate}
\item[$(a)$]  $\gT = \{\bar t:\bar t = \langle t_n:n < \omega \rangle,
2 \le t_n \in \bbN\}$,
\sn
\item[$(b)$]   we call $\bar t \in \gT$ full when $(\forall k \ge
  2)(\exists n)[k$ divide $\prod\limits_{\ell < n} t_\ell]$,
  equivalently $(\forall n)(\exists m)[m > n \wedge
  n|\prod\limits^{m}_{\ell = n} t_\ell]$, equivalently, every prime
  $p$, divide infinitely many $t_n$'s
\sn
\item[$(c)$]  we call $\bar t \in \gT$ explicitly weakly full 
\when \, for every
prime $p$, either $p$ divide no $t_n$ or it divides infinitely many
$t_n$
\sn
\item[$(d)$]  we say $G$ is $\bar t$-divisible \when \, every $x \in
G$ is divisible by $\prod\limits_{\ell < n} t_\ell$ for every $n$
\sn
\item[$(e)$]  we call $\bar t \in \gT$ weakly full \when \, for some
  $n(*)$ the sequence $\langle t_{n(*)+n}:n < \omega\rangle$ is
  explicitly weakly full.
\end{enumerate}
\mn
3) For $G \in K^{\rtf}_{\bar t,\lambda}$ let $G^{[\bar t]}$ 
be the $d_{\bar t}$-completion of $G$ \underline{where} 
$d_{\bar t} = d_{\bar t}[G]$
is the metric defined by $d_{\bar t}(x,y) = \inf\{2^{-k}:
\prod\limits_{\ell < k} t_\ell$ divides $x-y$ in the Abelian group 
$G\}$, justify by
\ref{a12}(3), pedantically ``the $d_t$-completion" is 
determined only up to isormophism over $G$.

\noindent
4) Let $K^{\crtf}_{\bar t,\lambda}$ be the class of $G \in 
K^{\rtf}_{\bar t,\lambda}$ which are $\bar t$-complete 
(i.e. $G = G^{[\bar t]}$).

\noindent
5) For those classes, $\le$ means being a subgroup and
$\le_{\pr}$ means being a pure subgroup.

\noindent
6) We say $\bar t,\bar s \in \gT$ are equivalent when
$K^{\rtf}_{\bar t} = K^{\rtf}_{\bar s}$.
\end{definition}

\begin{observation}
\label{a12}
1) $\bar t$ is full iff $\bar t$ is equivalent to $\langle n!:n \in
\bbN\rangle$ \Iff \, for every power of prime $m$, for some
$n,m$ divides $\prod\limits_{\ell < n} t_\ell$.

\noindent
2) If $\bar t$ is full \then \, every $G \in K^{\tf}$ can be
represented in fact uniquely as the direct sum 
$G_1 + G_2$ where $G_1$ is divisible, $G_2 \in K^{\rtf}_{\bar t}$.

\noindent
3) For $G \in K^{\rtf}_{\bar t},d_{\bar t}$ is a metric on $G$.

\noindent
4) If $G \in K^{\rtf}_{\bar t}$ \then \, there is $G'$, called the
$\bar t$-completion of $G$, such that
\mn
\begin{enumerate}
\item[$(a)$]  $G \le_{\pr} G' \in K^{\rtf}_{\bar t}$
\sn
\item[$(b)$]  $G'$ is $\bar t$-complete
\sn
\item[$(c)$]  $G$ is dense in $G'$ by the metric $d_{\bar t}$
\sn
\item[$(d)$]  if $G''$ satisfies (a),(b),(c) \then \, $G'',G'$ are
isomorphic over $G$.
\end{enumerate}
\mn
5) $\bar t,\bar s \in \gT$ are equivalent \when \, for some $k,\ell$
we have
\mn
\begin{enumerate}
\item[$\bullet$]  $t_{k+n} = t_{\ell +n}$ for every $n$ or at least
\sn
\item[$\bullet$]  for some $m_*$, for every $m \ge m_*$ there is $n$ such that
$\prod\limits^{m}_{i = m_*} t_{k+i}$ divide $\prod\limits_{i<n} s_{\ell +i}$
and $\prod\limits_{\ell < m} s_{\ell +i}$ divides $\prod\limits_{i<n}
t_{\ell +i}$. 
\end{enumerate}
\mn
6) For members of $\gT$ being full and being weakly full 
are preserved by equivalence.
\end{observation}
\bigskip

\begin{proof}
\underline{Proof of \ref{a12}}:

Should be clear.
\end{proof}

\begin{PROOF}{\ref{a6}}

\noindent
\underline{Proof of \ref{a6}}

Let $t_n = n!$ and let $\bar t = \langle t_n:n < \omega \rangle$.

\noindent
The point is that clearly
\mn
\begin{enumerate}
\item[$(a)$]   $(\alpha) \quad$ for 
$G \in K^{\rtf}_{\bar t},G \le_{\pr} G^{[\bar t]} 
\in K^{\rtf}_{\bar t}$ and $G^{[\bar t]}$ has cardinality 
$\le \|G\|^{\aleph_0}$

\hskip30pt  and $G^{[\bar t]}$ is $d_{\bar t}$-complete, 
remember $G^{[\bar t]}$ is the $d_{\bar t}$-completion of $G$,

\hskip30pt  it is unique up to isomorphism over $G$
\sn
\item[${{}}$]  $(\beta) \quad$ if $G_1 \le_{\pr} G_2$ then $G^{[\bar
t]}_1 \le_{\pr} G^{[\bar t]}_2$, more pedantically: if $G_1 \le_{\pr}
G_2$

\hskip30pt $\le_{\pr} G_3$ and $G_3$ is $\bar t$-complete \then \, 
$G^{[\bar t]}_1$ can be (purely) embedded 

\hskip30pt into $G_3$ over $G_1$.
\end{enumerate}
\mn
Recall $K^{\crtf}_{\bar t}$ is the class of $d_{\bar t}$-complete
$G \in K^{\rtf}_{\bar t}$.  

Easily:
\mn
\begin{enumerate}
\item[$(b)$]    $(K^{\crtf}_{\bar t},\le_{\pr})$ has amalgamation, the
  joint embedding property and the LST (= L\"owenheim-Skolem-Tarski) property 
down to $\lambda$ for any $\lambda = \lambda^{\aleph_0}$
\sn
\item[$(c)$]   if $G' \le_{\pr} G''$ are from $K^{\ctrf}$ 
then we can find $\le_{\pr}$-increasing
sequence $\langle G_\alpha:\alpha \le \alpha(*)\rangle$ of members of
$K^{\crtf}$ such that
\begin{enumerate}
\item[$(\alpha)$]   $G' = G_0,G'' = G_{\alpha(*)}$
\sn
\item[$(\beta)$]  $x_\alpha \in G_{\alpha +1} \backslash G_\alpha$
\sn
\item[$(\gamma)$]   $G_{\alpha +1}$ is the $\bar t$-completion 
of the pure closure of $G_\alpha \oplus \bbZ x_\alpha$ inside
$G_{\alpha +1}$
\sn
\item[$(\delta)$]   for $\alpha$ limit, $G_\alpha$ is the $\bar
t$-completion of $\cup\{G_\beta:\beta < \alpha\}$ inside $G''$, note
that if $\cf(\alpha) > \aleph_0$ then the union is $\bar t$-complete.
\end{enumerate}
\sn
\item[$(d)$]  if $\lambda = \lambda^{\aleph_0}$ then 
for each $G \in K^{\crtf}_{\bar t,\le \lambda}$, we can
find $\langle (G_i,x_j):i \le \lambda^{\aleph_0},j < 
\lambda^{\aleph_0} \rangle$ such that:
\sn
\begin{enumerate}
\item[$(\alpha)$]   $G_0 = G,G_i$ is $\le_{\pr}$-increasing
continuous,
\sn
\item[$(\beta)$]  $x_i \in G_{i+1} \in K^{\crtf}_{\bar t,\lambda}$
\sn
\item[$(\gamma)$]  letting $G'_i$ be the pure closure of $G + \bbZ
x_i$ inside $G_* = \cup\{G_j:j < \lambda^{\aleph_0}\}$, we have
$G_{i+1} = G_i \bigoplus\limits_{G} G'_i$
\sn
\item[$(\delta)$]   if $G \le_{\pr} G',x \in G' \in 
K^{\crtf}_{\bar t,\lambda}$ and $G'$ is the $\bar t$-completion
of the pure closure of $G + \bbZ x$ inside $G'$
\then \, we can find $i < \lambda^{\aleph_0}$ and a 
pure embedding $h$ of $G'$ into $G_{i+1},h \restriction G =$ 
the identity, $h(x) = x_i$ (so $h''(G_i) \le_{\pr} G$), in fact, $h$
is onto $G'_i$.
\end{enumerate}
\sn
\item[$(e)$]   if $\lambda,G$ are as in clause (d) then we can find 
$G_* = \cup\{G_i:i < \lambda^{\aleph_0}\}$ such that
\sn
\begin{enumerate}
\item[$(\alpha)$]  $G \le_{\pr} G_* \in K^{\rtf}_{\lambda^{\aleph_0}}$
\sn
\item[$(\beta)$]  if $G \le_{\pr} G' \in K^{\rtf}_{\lambda^{\aleph_0}}$
then $G'$ can be purely embedded into $G_*$ over $G$
\sn
\item[$(\gamma)$]
\begin{itemize}
\item   $\langle G_i:i < \lambda^{\aleph_0}\rangle$ is a
$\le_{\pr}$-increasing continuous sequence of members of 
$K^{\rtf}_{\lambda^{\aleph_0}}$
\sn
\item   $G_0 = G$
\sn
\item   for limit $\delta <
\lambda^{\aleph_0},G_{\delta +1}$ is the $\bar t$-completion of
$G_\delta$
\sn
\item   for non-limit $\alpha < \lambda^{\aleph_0}$, the pair
$(G_\alpha,G_{\alpha +1})$ is like $(G,G_*)$ in clause (d)
\end{itemize}
\end{enumerate}
\sn
\item[$(f)$]  if for $i=1,2,G_\ell \in K^{\ctrf}_{\bar t,\lambda}$ and $\langle
G^\ell_i:i < \lambda^{\aleph_0}\rangle,G^\ell_*$ are as
in clause (d) or as in clause (e) and $\pi$ is an isomorphism from
$G_1$ onto $G_2$ \then \, there is an isomorphism $\pi^+$ from $G^1_*$
onto $G^2_*$ extending $\pi$
\sn
\item[$(g)$]  if $\lambda = \Sigma\{\lambda_n:n < \omega\},\lambda_n
= \lambda^{\aleph_0}_n < \lambda_{n+1}$ and $G \in 
K^{\rtf}_{\le\lambda}$ then we can find $G',G'_n$ such that
\sn
\begin{enumerate}
\item[$(\alpha)$]    $G \le_{\pr} G' \in K^{\rtf}_\lambda$
\sn
\item[$(\beta)$]    $G'_n \in K^{\crtf}_{\lambda_n}$
\sn
\item[$(\gamma)$]    $G'_n \le_{\pr} G'_{n+1}$; moreover there is a sequence
$\langle G'_{n,i}:i < \lambda^{\aleph_0}_n\rangle$ as in (e)
for $G'_n$ such that $G'_{n+1} = \cup\{G'_{n,i}:i < \lambda^{\aleph_0}_n\}$
\sn
\item[$(\delta)$]    $G' = \cup\{G'_n:n < \omega\}$.
\end{enumerate}
\sn
\item[$(h)$]  give $\lambda,\lambda_n$ as in (g), if $G',G''$ are as
$G'$ is in (g) \then \, $G',G''$ are isomorphic
\sn
\item[$(i)$]  moreover if $\lambda,\lambda_n$ are as in clause (g) and
$H \in K^{\rtf}_{\le \lambda}$ \then \, $H$ can be purely embedded
into $G'$; (and if $H \supseteq G$ then even embedded over $G$).
\end{enumerate}
\mn
The results now follow. 
\end{PROOF}

\noindent
In \ref{a18}(2) below we prove that there is no universal in $\lambda =
\lambda^{\aleph_0}$, using \cite[Th.1.1]{Sh:309}, for the reader's
convenience we quote the special case used.
\begin{fact}
\label{a15}
For any $\lambda$ and $X$, a set of cardinality $\le \lambda$ or just
$\le \lambda^{\aleph_0}$ \then \, we can find a sequence $\bar f = \langle
f_\eta:\eta \in {}^\omega \lambda\rangle$ such that:
\mn
\begin{enumerate}
\item[$(a)$]  $f_\eta$ is a function from $\{\eta \rest n:n <
  \omega\}$ into $X$
\sn
\item[$(b)$]  if $f$ is a function from ${}^{\omega >}\lambda$ to $X$
  \then \, for some $\eta \in {}^\omega \lambda$ we have $f_\eta
  \subseteq f$.
\end{enumerate}
\end{fact}

\begin{remark}
1) Concerning \ref{a15}, see \cite[1.5]{Sh:309}.

\noindent
2) We use \ref{a15} mainly for $\lambda = \lambda^{\aleph_0}$.
\end{remark}

\begin{claim}
\label{a18}
Assume $\bar t \in \gT$ is not full.

\noindent
1) $(K^{\rtf}_{\bar t},\le_{\pr})$ fails amalgamation.

\noindent
2) If $\lambda = \lambda^{\aleph_0}$ \then \, in 
$(K^{\rtf}_{\bar t,\lambda},\le_{\pr})$ there is no universal 
member, even for the $\aleph_1$-free ones.
\end{claim}

\begin{remark}
Note that \ref{a6}, \ref{a18}(2) are not contradictory as the former
deals with full $\bar t$ and the latter wtih non-full ones.
\end{remark}

\begin{PROOF}{\ref{a18}}
Let $p$ be a prime witnessing $\bar t$ is not full, i.e. $n_*$ is well
defined where $n_*= \min\{n:p$ divide no $t_m$ with 
$m \ge n\}$, by \ref{a12}(5) \wilog \, $n_*=0$.

Let $t_{<n} := \prod\limits_{\ell < n} t_\ell$ so $t_{<0} = 1$.

We now choose $a^1_n,a^0_n$ by induction on $n$ such that
\mn
\begin{enumerate}
\item[$(*)_1$]
\begin{enumerate}
\item[(a)]  $a^1_n,a^0_n \in \bbZ$
\sn
\item[(b)]  $a^1_n = a^0_n \mod t_{\le n}$
\sn
\item[(c)]  $a^\ell_n = a^\ell_m \mod p t_{\le n}$ if $n =m+1$
\sn
\item[(d)]  $a^1_n \ne a^0_n \mod p$ if $n=0$.
\end{enumerate}
\end{enumerate}
\mn
[Why we can choose?  For $n=0$ clearly $t_{<0}=1$ 
hence $a^1_n = 1,a^0_n=t_0$ are as required.

For $n=m+1$ let $a^0_n = a^0_m,a^1_n = a^1_m - p t_{<m} b_n$ for $b_n$
chosen below.  So clause $(*)_1(c)$ holds for $\ell = 0$ trivially and
for $\ell=1$ obviously.  Also clause $(*)_1(b)$ means $(a^1_m - a^0_n)
= p t_{<n} b_n,\mod t_{\le n}$.  By the induction hypothesis $b'_m =
(a^1_m - a^0_n)/t_{<n} \in \bbZ$ so the $(*)_1(b)$ means $b'_n = p t_n
\mod t_n$; as $p|t_n$ there is such $b_m$.

Lastly, $(*)_7(d)$ holds obviously by $(*)_2(b)$ and $(*)_1(d)$ for
$n=0$.]

Choose
\mn
\begin{enumerate}
\item[$(*)_2$]  $(a) \quad t'_n$ is $p t_n$ if $n = 0$ and is $t_n$ if
  $n>0$
\sn
\item[${{}}$]  $(b) \quad t'_{<n} = \prod\limits_{k<n} t'_k$ and
  $t'_{\le n} = t'_{<(n+1)}$
\sn
\item[${{}}$]  $(c) \quad c^\ell_n \in \bbZ$ are chosen such that 
$\sum\limits_{m \le n} t'_{<m} c^\ell_m = a^\ell_n$.
\end{enumerate}
\mn
[Why we can choose?  Just choose $c^\ell_n$ by induction on $n$.

For $n=0$ let $c^\ell_n = a^\ell_n$.

For $n = m+1$ let $c^\ell_n = (a^\ell_n - a^\ell_m)/t'_{\le m}$ which
belongs to $\bbZ$ by $(*)_1(c)$, now check the equation \newline
$\sum\limits_{i \le n} t'_{<i} c^\ell_i = (\sum\limits_{i<n} t'_i
c^\ell_i) + t'_{<n} c^\ell_n = a^\ell_m + t'_{\le m} c^\ell_n$ which
by the choice above is equal to a $a^\ell_n$ as required.]

For every $S \subseteq {}^\omega \lambda$ we let $G_S$ the Abelian
group generated by $\{x_\eta:\eta \in {}^{\omega >}\lambda\} \cup
\{y_{\eta,n}:\eta \in {}^\omega \lambda$ and $n < \omega\}$ freely
except the equations:
\mn
\begin{enumerate}
\item[$(*)_3$]  $t'_n y_{\eta,n+1} = y_{\eta,n} - c^\ell_n
  x_{\langle \rangle} - x_{\eta \rest n}$ \If \, $n < \omega$ and
  $\eta \in S \Rightarrow \ell=1$ and $\eta \notin S \Rightarrow \ell=0$.
\end{enumerate}
\mn
Let
\mn
\begin{enumerate}
\item[$(*)_4$]  for $n \in {}^\omega \lambda$ let
\sn
\begin{enumerate}
\item[$(a)$]  $G_\eta = \Sigma\{\bbZ x_{\eta \rest n}:n < \omega\}
\subseteq G_S$
\sn
\item[$(b)$]  $G_{S,\eta} = \Sigma\{\bbZ x_{\eta \rest n}:n < \omega\}
  + \Sigma\{\bbZ y_{\eta,n}:n < \omega\} \subseteq G_S$.
\end{enumerate}
\end{enumerate}
\mn
Easily
\mn
\begin{enumerate}
\item[$(*)_5$]  if $S \subseteq {}^\omega \lambda$ then
\sn
\begin{enumerate}
\item[$(a)$]  $G_S \in K^{\rtf}_{\bar t,\lambda^{\aleph_0}}$ is 
$\aleph_1$-free
\sn
\item[$(b)$]   $\eta \in {}^\omega \lambda \Rightarrow G_\eta
  \le_{\pr} G_{S,n} \le_{\pr} G_S$.
\end{enumerate}
\end{enumerate}
\mn
Now
\mn
\begin{enumerate}
\item[$\boxplus$]  if $S_0,S_1 \subseteq {}^\omega \lambda,\eta \in
  S_1 \backslash S_0$ \then \, $G_{S_0},G_{S_1}$ and even
  $G_{S_0,\eta},G_{S_1,\eta}$ cannot be $\le_{\pr}$-amalgamated over $G_\eta$.
\end{enumerate}
\mn
Why?  Toward contradiction assume $G_\eta \le_{\pr} H \in K^{\rtf}_{\bar
    t}$ and $\pi_\ell$ is a pure embedding of $G_{S_\ell}$ into $H$
  over $G_\eta$, for $\ell = 0,1$.

Let $z_n = \pi_1(y_{\eta,n}) - \pi_0(y_{\eta,n})$ for any $n$ and let $\pi =
\pi_0 \rest G_\eta = \pi_1 \rest G_n$.  

For any $n$ clearly for $\ell = 1,2$
\mn
\begin{enumerate}
\item[$\bullet_1$]  $G_\ell \models t'_{\le n} y'_{\eta,n+1} =
  y_{\eta,0} - (\sum\limits_{m \le n} t'_{<m} c^\ell_m)
  x_{\langle \rangle} + \sum\limits_{m \le n} t'_{<m} x_{\eta \rest m}$.
\end{enumerate}
\mn
So applying $\pi_\ell$ on the equation recalling $(*)_2(c)$ we have
\mn
\begin{enumerate}
\item[$\bullet_2$]  $H \models \pi_\ell(t'_{\le n} y_{\eta,n+1}) =
\pi_\ell(y_{\eta,0}) - a^\ell_n \pi(x_{\langle \rangle}) -
\sum\limits_{m \le n} t'_{<m} \pi(x_{\eta \rest n})$.
\end{enumerate}
\mn
Subtracting the equation in $\bullet_2$ for $\ell=0,1$ 
recalling the choice of $z_0,z_n$
\mn
\begin{enumerate}
\item[$\bullet_3$]  $H \models t'_{\le n} z_{n+1} = z_0 - (a^1_n - a^0_n) 
\pi(x_{\langle \rangle})$.
\end{enumerate}
\mn
But $t'_{\le n}$ and $a^1_n - a^0_n$ are divisible by $t_{\le n}$ in $\bbZ$
(by $(*)_2(a),(b)$ and $(*)_1(c)$ respectively) hence
\mn
\begin{enumerate}
\item[$\bullet_4$]  $z_0$ is divisible by $t_{\le n}$ in $H$.
\end{enumerate}
\mn
As this holds for every $n$ and $H \in K^{\rtf}_{\bar t}$ we get
\mn
\begin{enumerate}
\item[$\bullet_5$]  $z_0 = 0_H$.
\end{enumerate}
\mn
So $H \models t'_{\le n} z_{n+1} = z_0 - 
(a^1_n-a^0_n) \pi(x_{\langle \rangle})$ and for
$n=1$ we get $H \models t'_0 z_1 = z_0 - (a^1_0-a^0_0) 
\pi(x_{\langle \rangle})$, but in $\bbZ$ we have $p | t'_0$ 
and $p \dag (a^1_0 -
a^1_1)$ and $z_0 = 0$ so $p$ divides $x_{\langle \rangle}$ in $H$,
contradiction to purity.

This is enough for part (1), for part (2) 
we apply the simple black box of \cite[Th.1.1]{Sh:309},
i.e. \ref{a15}.  In details assume $G_* \in K^{\rtf}_\lambda$ and let
$\bar f = \langle f_\eta:\eta \in {}^\omega \lambda\rangle$ be as in
Fact \ref{a15} for $X = G_*$.

Define $S$ as the set of $\eta \in {}^\omega \lambda$ such that:
\mn
\begin{enumerate}
\item[$\bullet$]   there is no pure embedding
$g_*$ of $G_{\{\eta\},\eta}$ into $G_*$ such that $n < \omega
\Rightarrow g_*(x_{\eta \rest m}) = f_\eta(\eta \rest n)$.
\end{enumerate}
\mn
Now $G_S \in K^{\rtf}_{\bar t,\lambda}$ so it is
enough to prove that $G_S$ is not purely embeddable into $G_*$.
Toward contradiction assume $g$ is a pure embedding of $G_S$ into
$G_*$ and let $f:{}^{\omega >}\lambda \rightarrow X = G_*$ be $f(\eta)
= g(x_\eta)$.  By the choice of $\bar f$ there is $\eta \in {}^\omega
\lambda$ such that $f_\eta \subseteq f$.  If $\eta \in S$ then
$g \rest G_{S,\eta} = g \rest G_{\{\eta\},\eta}$ 
witness that $\eta \notin S$ by the definition of $S$.

So necessarily $\eta \in {}^\omega \lambda \backslash S$, hence 
there is $g_*$ as
forbidden in the definition of $S$.  Let $g_0 = g \rest G_{S,\eta}$.
Easily this contradicts $\boxplus$.
\end{PROOF}

\begin{remark}
\label{a21}
1) See more in \cite[Ch.II,\S3]{Sh:300} = \cite{Sh:300b}.

\noindent
2)  This holds also for $K^{\rs(p)}_\lambda$ the class
of reduced separable Abelian $p$-groups see \ref{a40}. 
\end{remark}

\noindent
We may wonder: what if we ask about $(K^{\rtf}_{\bar t,\lambda},\le)$,
i.e. the embedding is not necessarily pure.
\begin{claim}
\label{a25}
Assume $\bar t \in \gT$ is weakly full so
for some $n_*$ we have: if a prime $p$ divises some $t_n,n \ge n_*$ then
it divides infinitely many $t_n$'s, call this set of primes $\mathbf P$.

\noindent
1) If $\lambda = \lambda^{\aleph_0}$ then $(K^{\rtf}_{\bar
  t,\lambda},\le)$ has a universal member.

\noindent
2) If $\lambda = \sum\limits_{n} \lambda_n,\lambda_n =
(\lambda_n)^{\aleph_0}$ for every $n$ \then \, $(K^{\rtf}_{\bar
  t,\lambda},\le)$ has a universal member.

\noindent
3) Let $R$ be the subring of $\bbQ$ generated by $\{1\} \cup \{1/p:p$
a prime $\notin \mathbf P\}$.  Then for every $G \in K^{\rtf}_{\bar
  t,\lambda}$ there is $H \in K^{\trf}_{\bar t,\lambda}$ extending $G$
which is $p$-divisible for every prime $p \notin \mathbf P$.  Hence $H$
can be considered an $R$-module.

\noindent
4) For the class of $R$-modules into which $\bbQ_R$ cannot be
embedded the results of 1),2) holds replacing $\aleph_0$ by $|R| +
\aleph_0$ \when \, $R$ is an integral domain which is not a field,
$\bbQ_R$, its ring of quotients.
\end{claim}

\begin{PROOF}{\ref{a25}}
1),2)  By 4) and 3).

\noindent
3) Easy.

\noindent
4) The proof is like the proof for full $\bar t$.  
\end{PROOF}

\noindent
Still leaves some $\bar t$'s open.
\begin{claim}
\label{a28}
Assume $\bar t \in \gT$ is not weakly full hence $\mathbf P
:= \{p:p$ a prime dividing some $t_n$'s but only finitely many $t_n$'s$\}$ is
infinite (this is the negation of the conditions from \ref{a25}).

\noindent
If $\lambda = \lambda^{\aleph_0}$ then $(K^{\rtf}_{\bar
  t,\lambda},\le)$ has no universal member.
\end{claim}

\begin{PROOF}{\ref{a28}}
By \ref{a12}(5) \wilog 
\mn
\begin{enumerate}
\item[$(*)_1$]  $(a) \quad$ there are distinct primes $p_n$ such that:
  $p_k \big | t_n$ \Iff \, $k=n$
\sn
\item[${{}}$]  $(b) \quad (p_k)^{\ell(k)}$ divide $t_k$ but
$(p_k)^{\ell(k)+1}$ does not, so $\ell(k) \ge 1$.
\end{enumerate}
\mn
Let $t_{<n} = \prod\limits_{\ell < n} t_\ell$, so $t_{<0}=1$ and let
$t'_n = t_n p^{\ell(n)}_n$ and $t'_{<n} = \prod\limits_{\ell < n}
t'_\ell$ and 
$t''_n = p^{\ell(n)}_n$ and $t''_{<n} = \prod\limits_{\ell < n} t''_\ell$;
let $(t_{\le n},t'_{\le n},t''_{\le n}) =
(t_{<(n+1)},t'_{<(n+1)},t''_{<(n+1)})$.  

We now choose
$a^1_n,a^0_n \in \bbZ$ by induction on $n$ such that:
\mn
\begin{enumerate}
\item[$(*)_2$]  $(a) \quad a^1_n,a^0_n \in \bbZ$
\sn
\item[${{}}$]  $(b) \quad a^1_n = a^0_n \mod t'_{< m}$
\sn
\item[${{}}$]  $(c) \quad a^\ell_n = a^\ell_m \mod t'_{\le m}$
\sn
\item[${{}}$]  $(d) \quad$ if $k<n$ then $a^1_n \ne a^0_n \mod
  (p_k)^{\ell(k)+1}$.
\end{enumerate}
\mn
[Why possible?  First, for $n=0$ let $(a^1_n,a^0_n) = (p_0,t_0,t_0)$
so $a^1_n - a^0_n$ is divisible by $t_0$ but not by
$p^{\ell(n)+1}_n$.  Second, assume $n=m+1$ and $(a^1_m,a^0_m)$ have
been chosen.  As $t_{\le m}(t_{\le n}$ and $k \le m \Rightarrow p_k
\pm (t_{\le n}/t_{\le m})$ we can find $(b^1_m,b^0_m)$ such that
$b^\ell_m = a^\ell_m \mod t^*_{\le m}$ for $\ell=0,1$ and $b^1_m =
b^0_m \mod t^*_{\le n}$.  Clearly requirement (a),(b),(c) holds and
(d) for $k < m$.  Let $(a^1_n,a^0_n) = (a^1_n + t^*_{\le m} \cdot
t_n,a^0_n)$, now check.]
\mn
\begin{enumerate}
\item[$(*)_3$]  choose $c^1_n,c^0_n$ by induction on $n$ such that for 
$\ell=0,1$ we have $\sum\limits_{m \le n} t'_{< m}
c^\ell_n = a^\ell_n$.
\end{enumerate}
\mn
[Why possible?  For $n=0$ trivial for $n+1$ note the $c^\ell_{n+1}$
appear with coefficient 1.]

Next for every $S \subseteq {}^{\omega >}\lambda$ we choose an Abelian
group $G_S$, it is generated by $\{x_\eta:\eta \in {}^{\omega
  >}\lambda\} \cup \{y_{\eta,n}:\eta \in {}^\omega \lambda$ and $n <
\omega\} \cup \{x^*_n:n < \omega\}$ freely except the equations
\mn
\begin{enumerate}
\item[$(*)_4$]  
\begin{enumerate}
\item[(a)]  $(t_n/p^{\ell(n)}_n) x^*_{n+1} = x^*_n$ 
and $x_{\langle \rangle} = x^*_0$
\sn
\item[(b)]  $t'_n y_{\eta,n+1} = y_{\eta,n} - c^\ell_n
  x_{\langle \rangle} + x_{\eta \rest n}$ when $n < \omega,\ell <
  2$ and $\ell=1 \Leftrightarrow \eta \in S$ .
\end{enumerate}
\end{enumerate}
\mn
Also
\mn
\begin{enumerate}
\item[$(*)_5$]  
\begin{enumerate}
\item[(a)]   for $\eta \in {}^\omega \lambda$ let 
$G_\eta = \sum\limits_{n} \bbZ x_{\eta \rest n} + \sum\limits_{n} \bbZ
x^*_n$
\sn
\item[(b)]  for $S \subseteq {}^\omega \lambda,\eta \in
  {}^\omega \lambda$ let $G_{S,\eta}$ be the following subgroup of
$G_S:G_\eta + \Sigma \bbZ y_{\eta \rest n}$.
\end{enumerate}
\end{enumerate}
\mn
Easily
\mn
\begin{enumerate}
\item[$(*)_6$] 
\begin{enumerate}
\item[(a)]  if $S \subseteq {}^\omega \lambda$ and 
$\eta \in {}^\omega \lambda$ then 
$G_\eta,G_S,G_{S,n} \in K^{\rtf}_{\bar t,\lambda^{\aleph_0}}$
\sn
\item[(b)]  $G_\eta \le_{\pr} G_{S,\eta} \le_{\pr} G_S$.
\end{enumerate}
\end{enumerate}
\mn
[Why?  Note that $\bigcup\limits_{n} \bbZ_n x^*_n \in K^{\rtf}_{\bar
  t}$ because for every $n,x^*_0 = x_{\langle \rangle}$ is not
divisible by $p_n$.]

Now
\mn
\begin{enumerate}
\item[$\boxplus$]  if $S_0,S_1 \subseteq {}^\omega \lambda$ and $\eta
  \in S_1 \backslash S_0$ \then \, $G_{S_0},G_{S_1}$ and even
  $G_{S_0,\eta},G_{S_1,\eta}$ cannot be amalgamated over $G_\eta$ in
  $(K^{\rtf}_{\bar t_*},\le)$.
\end{enumerate}
\mn
We continue as in the proof of \ref{a18}, getting
$\pi_1,\pi_0,\pi,\eta,z_n$ and proving that for every $n$
\mn
\begin{enumerate}
\item[$\bullet$]  $H \models t'_{\le n} z_{n+1} = z_0 - (a^1_n - a^0_n) 
\pi(x_{\langle \rangle})$.
\end{enumerate}
\mn
But $t'_{\le n}$ is divisible by $t_{\le n}$ and $(a^1_n - a^0_n)$ is divisible
by $t''_{\le n}$ (in $\bbZ$) and $x_{\langle \rangle}$ is divisible 
by $t_{\le n}/t''_{\le n}$ hence
$(a^1_n-a^0_n) x_{\langle \rangle}$ is divisible by $t_{\le n}$, hence
$(a^1_n - a^0_n) \pi(x_{\langle \rangle})$ is divisible by $t_{\le
  n}$.  So by the equation above
$z_0 \in t_{<n} H$ for every $n$.  As $H \in K^{\rtf}_{\bar t}$ it
follows that $z_0 = 0$.

Hence for every $n$
\mn
\begin{enumerate}
\item[$\bullet$]  $H \models (a^1_n - a^0_n) \pi(x_{\langle \rangle}) 
= -t'_{\le n} z_n$.
\end{enumerate}
\mn
Now $p^{\ell(n)+\ell(n)}_n$ divide $t'_{\le n}$ and $p^{\ell(n)+1}_n$
does not divide $(a^1_{n+1} - a^0_{n+1})$ by $(*)_2(d)$ hence in
$H,p^{\ell(n)}_n$ divide $\pi(x_{\langle \rangle})$.  As also each 
$t'_{\le n}/\prod\limits_{k \le n} p^{\ell(k)}_k$ divide it, clearly
$\pi(x_{\langle \rangle})$ contradict $G_* \in K^{\rtf}_{\bar t,\lambda}$.
\end{PROOF}

\noindent
We may wonder whether the existence result of \ref{a6} holds for a
stronger embeddability notion.  A natural candidate is
\begin{definition}
\label{a34}
Let $G_0 \le_{\bar t} G_1$ if: $G_0,G_1$
are Abelian groups on which $\|-\|_{\bar t}$ is a norm, $G_0 
\le_{\pr} G_1$ and $G_0$ is a $d_{\bar t}$-closed subset of 
$G_1$ (but $G_\ell$ is not necessarily $\bar t$-complete!). 
\end{definition}

\begin{observation}
\label{a35}
1) $(K^{\rtf}_t,\le_{\bar t})$ is an a.e.c. except smoothness with LST
number $2^{\aleph_0}$.

\noindent
2) If $A \subseteq G \in K^{\rtf}_{\bar t}$ then for some $G'
\le_{\bar t} G,A \subseteq G',|G'| = (|A| + \aleph_0)^{\aleph_0}$.

\noindent
3) If $G_1 \le_{\bar t} G_2$ then $G_1 \le_{\pr} G_2$.
\end{observation}

\noindent
We prove below that for $\mu$ strong limit of cofinality $\aleph_0$
the answer is positive, i.e. there is a universal member for
$(K^{\rtf}_{\bar t,\lambda},\le_{\bar t})$, but for cardinals like 
$\beth^+_\omega < (\beth_\omega)^{\aleph_0}$ the question on the
existence of universals remain open.
\begin{fact}
\label{a37}
Assume $\lambda$ is strong limit and $\aleph_0 = \cf(\lambda) < 
\lambda$.

\noindent
1) There is a universal member in $(K^{\rtf}_{\bar t,\lambda},<_{\bar
  t})$ where $\bar t = \langle t_\ell:\ell < \omega \rangle 
\in \gT$, hence also in $(K^{\rtf}_{\bar t,\lambda},\le_{\pr})$ and
$(K^{\rtf}_{\bar t,\lambda},\subseteq)$.

\noindent
2) For a prime number $p$, similarly for 
$(K^{\rs(p)}_\lambda,\le_{\langle p:\ell <
   \omega\rangle})$, see Definition \ref{a40} below.
\end{fact}

\begin{definition}
\label{a40}
For a prime number $p$, and cardinal $\lambda$ we let
$K^{\rs(p)}_\lambda$ be the class of Abelian $p$-groups which 
are reduced and separable of cardinality $\lambda$.
\end{definition}

\begin{PROOF}{\ref{a37}}
Let $K$ be the class and $\le_*$ the partial order.
Let $\lambda_n < \lambda_{n+1} 
< \lambda = \sum\limits_{n} \lambda_n$ and $2^{\lambda_n} < 
\lambda_{n+1}$.  The idea in both cases 
is to analyze $M \in K_\lambda$ as the union of increasing chain 
$\langle M_n:n < \omega \rangle,M_n 
\prec_{\bbL_{\lambda^+_n,\lambda^+_n}} M,\|M_n\| = 2^{\lambda_n}$.

Specifically, we shall apply \ref{a43}, \ref{a46} below with:

\[
{\gK} = K^{\rtf}
\]

\[
\mu_n = (2^{\lambda_n})^+
\]

\[
\le_1 \text{ is}: M_1 \le_1 M_2 \text{ \underline{iff} }  (M_1,M_2 \in
K \text{ and) } M_1 \le_* M_2
\]

\begin{equation*}
\begin{array}{clcr}
\le_2 \text{ is}: M_1 \le_2 M_2 \text{ \underline{iff} } &M_1 \le_1
M_2 \text{ and}\\
  &M_1 \prec_{\bbL_{\aleph_1,\aleph_1}} M_2, \text{ or just}: \\
  &\text{if } G_1 \subseteq M_1,G_1 \subseteq G_2 \subseteq M_2, \\
  &\text{and } G_2 \text{ is countable then there is an} \\
  & \le_1\text{-embedding } h \text{ of } G_2 \text{ into } M_1 \text{ over } G_1.
\end{array}
\end{equation*}

\mn
We should check the conditions in \ref{a43} which we postpone.

\noindent
We shall finish the proof after \ref{a46} below.
\end{PROOF}

\begin{claim}
\label{a43}
Assume
\mn
\begin{enumerate}
\item[$(a)$]   $K$ is a class of models of a fixed vocabulary
  closed under isomorphism, $K_\lambda \ne \emptyset$
\sn
\item[$(b)$]   $\lambda = \sum\limits_{n < \omega} \mu_n,\mu_n <
\mu_{n+1},2^{\mu_n} < \mu_{n+1},\mu_n$ is regular and the 
vocabulary of $K$ has cardinality $< \mu_0$.
\sn
\item[$(c)$]   $\le_1$ is a partial order on $K$, (so $M \le_1 M$)
preserved under isomorphisms, and if $\langle M_i:i < \delta \rangle$ is
$\le_1$-increasing and continuous then $M_\delta = 
\bigcup\limits_{i < \delta} M_i \in K$ and $i < \delta 
\Rightarrow M_i \le_1 M_\delta$ (so $(K,\le_1)$ satisfies a quite
weak version of a.e.c. see \cite{Sh:88r} = \cite{Sh:88})
\sn
\item[$(d)$]
\sn
\begin{enumerate}
\item[$(\alpha)$]  $\le_2$ is a two-place relation on 
$K$, preserved under isomorphisms
\sn
\item[$(\beta)$]   [weak LST] if $M \in K_\lambda$ then 
we can find $\langle M_n:n < \omega \rangle$ such that:
$M_n \in K_{< \mu_n},M_n <_2 M_{n+1}$ and $M = 
\bigcup\limits_{M < \omega} M_n$
\end{enumerate}
\sn
\item[(e)]  [non-symmetric amalgamation]
 if $M_0 \in K_{< \mu_n},M_0 \le_1 M_1 \in K_{< \mu_{n +2}},
N^1 \le_2 N^2 \in K_{< \mu_{n+1}},h^1$ an
isomorphism from $M_0$ onto $N^1$, \then \, we can find $M_2 \in 
K_{< \mu_{(n+2)}}$ such that $M_1 \le_1 M_2$ and there is an 
embedding $h^2$ of $N^2$ into $M$ extending $h^1$ satisfying $h(N^2)
\le_1 M_2$.
\end{enumerate}
\mn
\Then \, we can find $\langle M^\alpha_n:n \le \omega \rangle$ for
$\alpha < 2^{< \mu_0}$ such that:
\mn
\begin{enumerate}
\item[$(\alpha)$]   $M^\alpha_n \in K_{< \mu_n},M^\alpha_n \le_1
M^\alpha_{n+1},M^\alpha_\omega = \bigcup\limits_{n < \omega} M^\alpha_n$
\sn
\item[$(\beta)$]   if $M \in K_\lambda$ and the sequence
$\langle M_n:n < \omega \rangle$ is as in clause $(d)(\beta)$
\then\footnote{Of course, we can omit $\langle M^\alpha_n:n \le
  \omega\rangle$ when $\|M^\alpha_\omega\| < \lambda$.} \, for 
some $\alpha < 2^{< \mu_0}$ we can find an embedding $h$ of $M$ 
into $M^\alpha_\omega$ satisfying $h(M_n) \le_1 M^\alpha_{n+2}$ (if
${\gK} = (K,\le_1)$ is an a.e.c. we get that $h$ is a $\le_{\gK}$-embedding 
of $M$ into $M^\alpha_\omega$).
\end{enumerate}
\end{claim}

\begin{PROOF}{\ref{a43}}
Let

\begin{equation*}
\begin{array}{clcr}
K'_0 = \bigl\{ M:&M \in K \text{ has universe an ordinal}\\
  &< \mu_0, \text{ and there is } \langle M_n:n < \omega \rangle \text{ as
in clause } (d)(\beta) \\
  &\text{with } M_0 \cong M \bigr\}.
\end{array}
\end{equation*}

\mn
Clearly $K'_0$ has cardinality $\le 2^{< \mu_0}$, and let us list it as
$\langle M^\alpha_0:\alpha < \alpha^* \rangle$ with $\alpha^* \le
2^{< \mu_0}$.  We now choose, for each $\alpha < \alpha^*$, by induction
on $n < \omega,M^\alpha_n$ such that:
\mn
\begin{enumerate}
\item[$(i)$]   $M^\alpha_n \in {\gK}$ has universe an ordinal $< \mu_n$
\sn
\item[$(ii)$]   $M^\alpha_n \le_1 M^\alpha_{n+1}$
\sn
\item[$(iii)$]   if $N^1 \le_2 N^2,N^1 \in K_{< \mu_n},N^2 \in
K_{< \mu_{n+1}},h^1$ is an embedding of $N^1$ into
$M^\alpha_{n+1}$ satisfying $h^1(N^1) \le_1 M^\alpha_{n+1}$ \then \,
we can find $h^2$, an embedding of $N^2$
into $M^\alpha_{n+2}$ extending $h^1$ such that $h^2(N^2) \le_1
M^\alpha_{n+2}$.
\end{enumerate}
\mn
For $n=0,1$ we do not have much to do.  (If $n=0$ use 
$M^\alpha_0$; if $n=1$ let 
$\langle M_n:n < \omega \rangle$ be as in clause (c), $M_0 \cong M^\alpha_0$
and use $M^\alpha_1$ such that $(M_1,M_0) \cong (M^\alpha_1,M^\alpha_0)$).
Assume $M^\alpha_{n+1}$ has been defined, and we shall define
$M^\alpha_{n+2}$, let $\{(h^1_{n,\zeta},N^1_{n,\zeta},N^2_{n,\zeta}):\zeta <
\zeta^*_n\}$ list the cases of clause
(iii) that need to be taken care of, with the set of elements of 
$N^2_{n,\zeta}$ being an ordinal.  \Wilog \, $\zeta^*_n \le 2^{<
  \mu_{n+1}}$ by cardinality consideration.  We shall choose 
$\langle N_{n+1,\zeta}:
\zeta \le \zeta^*_n \rangle$ which is $\le_1$-increasing continuous satisfying
$N_{n+1,\zeta} \in K_{< \mu_{n+2}}$.  We choose $N_{n+1,\zeta}$
by induction on $\zeta$.   Let
$N_{n+1,0} = M^\alpha_{n+1}$, for $\zeta$ limit let $N_{n+1,\zeta} =
\bigcup\limits_{\xi < \zeta} N_{n+1,\xi}$ and use clause (c) of 
the assumption.  

Lastly, for $\zeta = \xi +1$ use clause (e) of 
the assumption with $h^1_{n,\zeta}(N^1_{n,\xi}),N_{n+1,\xi},
N^1_{n,\xi},N^2_{n,\xi}$,
\newline
$h^1_{n,\xi},N_{n+1,\xi+1}$ 
here standing for $M_0,M_1,N^1,N^2,h^1,h^2,M_2$ there. 

Having carried the induction on $\zeta \le \zeta^*_n$ we let
$M^\alpha_{n+2} = N_{n+1,\zeta^*_\alpha}$; so we have carried the
induction on $n$.

Having chosen the $\left<\langle M^\alpha_n:n < \omega \rangle:\alpha
< 2^{< \mu_0}\right>$ let $M^\alpha_\omega = \cup\{M^\alpha_n:n <
\omega\}$ hence by clause (c) of the assumption, $M^\alpha_\omega \in
K$ and $n <\omega \Rightarrow M^\alpha_n \le_1
M^\alpha_\omega$.  Clearly clause $(\alpha)$ of the desired conclusion
is satisfied.  For clause $(\beta)$ let $M \in K_\lambda$.  By
clause (d) of the assumption we can find a sequence $\langle M_n:n <
\omega\rangle$ such that $M_n \in K_{< \mu_n},M_n \le_2
M_{n+1}$ and $M = \cup\{M_n:n < \omega\}$.  By the choice of $\langle
M^\alpha_0:\alpha < 2^{< \mu_0}\rangle$ there is $\alpha < 2^{<
\mu_0}$ such that $M_0 \cong M^\alpha_0$, and let $h_0$ be an
isomorphism from $M_0$ onto $M^\alpha_0$.  Now by induction on $n
<\omega$ we choose $h_n$, an embedding of $M_n$ into $M^\alpha_{n+1}$
such that $h_n(M_n) \le_1 M^\alpha_{n+1}$ and $h_n \subseteq
h_{n+1}$.  For $n=0$ this has already been done as $h_0(M_0) =
M^\alpha_0 \le_1 M^\alpha_1$.  For $n+1$ we use clause (iii).  

Lastly, $h =\cup\{h_n:n < \omega\}$ is an embedding of $M$ into
$M^\alpha_\omega$ as required.  
\end{PROOF}

\begin{remark}
\label{4.4A}
1) We can choose $\langle M^\alpha_0:\alpha < \alpha^* \rangle$ just 
to represent ${\gK}_{< \mu_0}$, 
and similarly later (and so ignore the ``with the universe being 
an ordinal"). 

\noindent
2) Actually, the family of $\langle M_n:n < \omega \rangle$ as in clause (c)
such that $M_n$ has set of elements an ordinal, forms a tree $T$ 
with $\omega$ levels with the 
$n$-th level having $\le 2^{< \mu_n}$ members, and we can use some 
amalgamations of it (so weakening the assumptions on $\le_1$).  
This gives a variant of \ref{a43}. 

\noindent
3) We can put into the axiomatization the stronger version of (d) 
from \ref{a43} proved in the proof of \ref{a37} so we can weaken 
$(\beta)$ of \ref{a46} below.

\noindent
4) E.g., in (d) we can add $M_n <_* M$ and so weaken clause $(\beta)$ 
of \ref{a43}.
\end{remark}

\begin{conclusion}
\label{a46}
1) In \ref{a43} we can add $\bigwedge\limits_{n} \,
   \bigwedge\limits_{\alpha} [M^\alpha_n = M^0_n]$ provided that:
\mn
\begin{enumerate}
\item[(f)]  there is $M_* \in K_{< \lambda}$ such that every $M
  \in K_{< \mu_0}$ can be $\le_1$-embeddable into $M_*$.
\end{enumerate}
\mn
2) In \ref{a43} there is in $K_\lambda$ a
universal member under $\le_1$-embedding \If \, in addition we add to
the assumptions of \ref{a43}: 
\mn
\begin{enumerate}
\item[$(f)^+$]   as in part (1)
\sn
\item[$(g)$]   if $M_n \le_1 M_{n+1},N_n \le_1 M_n,
N_n \le_2 N_{n+1}$ and $M_n \in K_{< \mu_{n+2}}$ and $N_n \in K_{<
\mu_{n+1}}$ for $n < \omega$ \then \, $\bigcup\limits_{n < \omega} N_n \le_1 
\bigcup\limits_{n < \omega} M_n$.
\end{enumerate}
\end{conclusion}

\begin{PROOF}{\ref{a46}}
Easy.
\end{PROOF}

\begin{remark}
\label{}
1) In \ref{a46}(2) we can weaken clause (f) to:
\mn
\begin{enumerate}
\item[(f)$'$]  there is $M_* \in K_{< \lambda}$ as there.
\end{enumerate}
\mn
2) This holds for \ref{a46}(1) as in \ref{a43} in clause $(\beta)$ we
replace $M^\alpha_{n+2}$ by $M^\alpha_{n+k}$ when $\|M_*\| < \mu_k$.
\end{remark}
\bigskip

\noindent
\underline{Continuation of the proof of \ref{a37}} 

We have to check the demands in \ref{a46} and \ref{a43}.

The least trivial clause to check is clause (e).
\bigskip

\noindent
\underline{Clause $(e)$}:  (non-symmetric amalgamation) 

Without loss of generality $h_1 =$ the identity, $N^1 \cap M_1 = M_0 = N_0$.
Just take the free amalgamation $M = N^1 *_{M_0} M_1$ (in the variety of
Abelian groups) and note that naturally $M_1 \le_1 M$.

\bigskip
\centerline {$* \qquad * \qquad *$}
\bigskip

\begin{discussion}
\label{a49}
1) Can we in \ref{a43}, \ref{a46} replace $\cf(\lambda) = 
\aleph_0$, by $\cf(\lambda) = \theta > \aleph_0$?  If
increasing union of chains in $K_{< \lambda}$ of length $< \theta$ 
behaves nicely then yes, with no real problem.  

More elaborately
\mn
\begin{enumerate}
\item[$(i)$]   in \ref{a43}(c), we get $\langle
M_\varepsilon:\varepsilon < \theta \rangle$ such that $M_\varepsilon
\in K_{< \mu_\varepsilon},\langle M_\varepsilon:\varepsilon < \theta
\rangle$ is $\subseteq$-increasing continuous, $M_\varepsilon <_2
M_{\varepsilon +1},M = \cup\{M_\varepsilon:\varepsilon < \theta
\rangle$
\sn
\item[$(ii)$]   we add: if $\langle M_i:i \le \delta \rangle$ is
$\le_1$-increasing continuous, $M_i \in K_{< \lambda}$ and $i < \delta
\Rightarrow M_i \le_1 N$ then $M_\delta \le_i N$.
\end{enumerate}
\mn
Otherwise we seem to be lost.  

\noindent
2) Suppose $\lambda = \sum\limits_{n < \omega} \lambda_n,\lambda_n =
(\lambda_n)^{\aleph_0} < \lambda_{n+1}$, and $\mu < \lambda_0,\lambda <
2^\mu$ (i.e., called Case 6b in \cite[\S0]{Sh:622}).  For $\bar t \in \gT$
which is not weakly full, is there a universal member in
$({\gK}^{\rtf}_{\bar t,\lambda},<_{\bar t})$? 

Assume $\mathbf V \models ``\mu = \mu^{<\mu},\mu < \chi"$ and 
$\bbP$ is the forcing notion of adding $\chi$ Cohen subsets 
to $\mu$ (that is $\bbP = \{f:f \text{ a partial function from } 
\chi \text{ to } 2, |\text{Dom}(f)| 
< \mu\}$ ordered by inclusion).  So we have in 
$\mathbf V^{\bbP}:\lambda < \lambda^{\aleph_0}$ and 
$\mu < \lambda < \chi \Rightarrow$ in $(K^{\rtf}_{\bar t,\lambda},
\le_{\bar t})$ there is no universal member.  The proof is easy; so
consistently the answer is no.

Maybe continuing \cite[\S2]{Sh:E59} =
\cite[Ch.III,\S2]{Sh:e} we can get consistency of the existence.

\noindent
3) Now if $\lambda = \lambda^{\aleph_0}$ then in 
$(K^{\aleph_1\text{-free}}_{\lambda},\subseteq)$ 
there is no universal member; see \cite{Sh:309} = \cite[Ch.IV]{Sh:e},
\cite{Sh:622} because amalgamation fails badly.  
Putting together those results clearly there are few cardinals which are
candidates for consistency of existence.  In (2), if there is a regular
$\lambda' \in (\mu,\lambda)$ with cov$(\lambda,\lambda^+,\lambda^+,\lambda')
< 2^\mu$ then contradict \ref{a6}.  

\noindent
4) Considering consistency of existence of universal in (2), it is 
natural to try to combine the
independent results in \cite{Sh:309} = \cite[Ch.IV]{Sh:e} and \cite{Sh:614}.
\end{discussion}
\newpage

\section {The class of Groups} \label{theclass}

We know (\cite{Sh:789}) that the class of groups has NSOP$_4$ and
SOP$_3$ (from \cite[\S2]{Sh:500}).  We shall prove a result on the
place of the class of groups in the model theoretic classification.
We know that it falls on ``the complicated side" for some division:
of course is unstable.  Now we prove that it has the 
oak property (see on it \cite{Sh:710}).
This is formally not well defined as the definition there was for 
complete first order theories.  But its meaning
(and ``no universal" consequences) are clear in a more general
context, see below.
Amenability is a condition on a theory (or class) which gives
sufficient condition for existence of somewhat
universal structures and in suitable
models of set theory (see \cite{Sh:614}), the class of groups fail
it because by \cite{Sh:1029} essentially, it has no universal in $\lambda$ when
$\lambda = \mu^+,\mu = \mu^{<\mu}$, forcing contradiction the results
on amenable elementary classes in \cite{Sh:614}.

\begin{definition}
\label{b3}
1)  A theory $T$ is said to satisfy the oak property
as exhibited by (or just by) a formula $\varphi(\bar x,\bar y,\bar z)$
\when \, for any $\lambda,\kappa$ there are $\bar b_\eta(\eta \in
{}^{\kappa >} \lambda)$ and $\bar c_\nu(\nu \in {}^\kappa \lambda)$ 
and $\bar a_i (i < \kappa)$ in some model $\gC$ of $T$ such that
\mn
\begin{enumerate}
\item[$(a)$]   $\eta \triangleleft \nu$ and $\nu \in {}^\kappa
  \lambda$ then $\gC \models \varphi[\bar a_{\ell g(\eta)},
\bar b_\eta,\bar c_\nu]$
\sn
\item[$(b)$]   if $\eta \in {}^{\kappa >} \lambda$ and $\eta \char 94
\langle \alpha \rangle \trianglelefteq \nu_1 \in {}^\kappa \lambda$ and $\eta
\char 94 \langle \beta \rangle \trianglelefteq \nu_2 
\in {}^\kappa \lambda$, while $\alpha
\ne \beta$ and $i > \ell g(\eta)$, \then \, $\neg \exists \bar
y[\varphi(\bar a_i,\bar y,\bar c_{\nu_1}) \wedge \varphi(\bar a_i,\bar
y,\bar c_{\nu_2})]$ 
\end{enumerate}
\mn
and in addition $\varphi$ satisfies
\mn
\begin{enumerate}
\item[$(c)$]   $\varphi(\bar x,\bar y_1,\bar z) \wedge
\varphi(\bar x,\bar y_2,\bar z)$ implies $\bar y_1 = \bar y_2$ in any
model of $T$.
\end{enumerate}
\mn
2) A theory $T$ has the $\Delta$-oak property if it is
exhibited by some $\varphi(\bar x,\bar y,\bar z) \in \Delta$.
\end{definition}

\begin{claim}
\label{b6}
The class of groups has the oak property by some quantifier free formula.
\end{claim}

\begin{remark}
The original proof goes as follows.

Let $w(x,y)$ be a complicated enough word, say of length $k^* = 100$,
see demands below.

\noindent
For cardinals $\kappa,\lambda$ let $G = G_{\lambda,\kappa}$ be defined
as follows:

Let $G$ be the group generated by $\{x_i:i < \kappa\} \cup
\{y_\eta:\eta \in {}^{\kappa >} \lambda\} \cup \{z_\nu:\nu \in
{}^\kappa \lambda\}$ freely except the set of equations

\[
\Gamma = \{y_{\nu \restriction i} = w(z_\nu,x_i):\nu \in {}^\kappa
\lambda,i < \kappa\}.
\]

\mn
Clearly it suffices to show that
\mn
\begin{enumerate}
\item[$(*)_1$]   if $\nu \in {}^\kappa \lambda,i < \kappa$ and $\rho \in
{}^i \lambda \backslash \{\nu \restriction i\}$ then $G \models ``y_\rho \ne
w(z_\nu,x_i)"$.
\end{enumerate}
\mn
Now
\mn
\begin{enumerate}
\item[$(*)_2$]   each word $y^{-1}_{\nu \restriction i}
w(z_\nu,x_i)$ is so-called cyclically reduced, i.e. 
both $w_1 = y^{-1}_{\nu \rest i}
w(z_v,x_i)$ and $w_2 = w(z_v,x_i)y^{-1}_{\nu \rest i}$ are reduced,
i.e. we do not have a generator and its inverse in adjacent places
\sn
\item[$(**)$]   for any two such words or cyclical permutations of
them which are not equal, any common segment has length $< k^*/6$.
\end{enumerate}
\mn
Explanation and why this is enough see \cite{LySc77}, no point to
elaborate as this is not used.

But we prefer to use the more ad-hoc but accessible proof.
\end{remark}

\begin{PROOF}{\ref{b6}}

\noindent
\underline{Proof of \ref{b6}}
Let $G = G_0$ be the group generated by

\[
Y = \{x_i:i < \kappa\} \cup \{z_\nu:\nu \in {}^\kappa \mu\}
\]

\mn
freely except (recalling $[xy] = xyx^{-1}y^{-1}$, the commutator)
the set of equations $\Gamma_2 = \{[z_\nu,x_i] =
[z_\eta,x_i]:i < \kappa,\nu \in {}^\kappa \lambda,\eta \in 
{}^\kappa \lambda$ satisfy $\nu \restriction i = \eta \restriction i\}$.  
So for $i < \kappa,\rho \in {}^i \lambda$ we
can choose $y_\rho \in G$ such that $\eta \in {}^\kappa \lambda,\eta
\restriction i = \rho \Rightarrow y_\rho = [z_\eta,x_i]$.  
Let $G_1$ be the group generated by the set $Y$ freely, let $h$ be the
homomorphism from $G_1$ onto $G$ mapping the members of $Y$ to
themselves, (using Abelian groups no two members of $Y$ are identified
in $G_1$).  Let $N =  \Kernel(h)$.

Clearly it suffices to prove
\mn
\begin{enumerate}
\item[$(*)_1$]   in $G=G_1/N$, if $\nu,\eta \in {}^\kappa \lambda$ and $i
< \kappa$ then $[z_\nu,x_i] = [z_\eta,x_i] \Leftrightarrow \nu
\restriction i = \eta \restriction i$.
\end{enumerate}
\mn
The implication $\Leftarrow$ holds trivially.  For the other direction
let $j<\kappa$ and $\eta,\nu \in {}^\kappa \lambda$ be such that $\eta
\restriction j \ne \nu \restriction j$ and we shall prove that $G
\models ``y_{\eta \restriction j} \ne y_{\nu \restriction j}"$. 

Let $N_1$ be the normal subgroup of $G_1$ generated by

\begin{equation*}
\begin{array}{clcr}
(*)_2 \qquad X_* = \{x_i:i < \kappa \text{ and } i \ne j\} &\cup \{z_\rho:\rho 
\in {}^\kappa \lambda \text{ and } \rho \restriction j \notin \{\eta 
\restriction j,\nu \restriction j\}\} \\
  &\cup \{z_\rho z^{-1}_\eta:\rho \in {}^\kappa \lambda \text{ and } \rho
  \restriction j = \eta \restriction j\} \\
  &\cup \{z_\rho z^{-1}_\nu:\rho \in {}^\kappa \lambda \text{ and } \rho
  \restriction j = \nu \restriction j\}.
\end{array}
\end{equation*}

\mn
Clearly by inspection $N_1$ includes $N$.  Let $N_0 = h(N_1)$, 
clearly $N_1$ is a
normal subgroup of $G_1$ and $h$ induces a homomorphism $\hat h$ from
$G_1/N_1$ onto $G_0/N_0$.  Now looking at the members of 
$X_*,G_1/N_1$ is generated by $\{x_i\} \cup
\{z_\eta,z_\nu\}$.  Checking the equations in $\Gamma_2$ clearly $G_1/N_1$ is
generated by $\{x_i\} \cup \{z_\eta,z_\nu\}$ freely, hence 
$G_1/N_1 \models ``[z_\eta,x_i] \ne [z_\nu,x_i]"$ which means
$[z_\eta,x_i]^{-1}[z_\nu,x_i] \notin N_1$ hence $\notin N$.  So recalling
the choice of $G$ in $(*)_1$ we have 
$G \models ``y_{\eta \restriction j} \ne y_{\nu \restriction j}"$ as required.
\end{PROOF}
\newpage

\section {More On the oak property} \label{ontheoak}

We can in the ``no universal" results in \cite{Sh:710} deal also
with the case of singular cardinals.  We also note that the so called
weak oak property suffices.
\begin{claim}
\label{c3}
We have $\univ(\lambda_1,T) \ge \lambda_2$ \when \,:
\mn
\begin{enumerate}
\item[$(a)$]   $T$ is a complete first order theory with the oak
property, ${\gK} = (\Mod_T,\prec)$ 
\sn
\item[$(b)$]
\begin{enumerate}
\item[$(i)$]  $\kappa = \cf(\mu) \le \sigma < \mu <
\lambda = \cf(\lambda) < \lambda_1 \le \lambda_2$
\sn
\item[$(ii)$]  $\kappa \le \sigma \le \lambda_1,|T| \le \lambda_2$
\sn
\item[$(iii)$]  $\cf([\mu]^\kappa,\subseteq) \ge \lambda_2$
\end{enumerate}
\sn
\item[$(c)$]
\begin{enumerate}
\item[$(i)$]  $S \subseteq \lambda$ is stationary
\sn
\item[$(ii)$]  $\bar C = \langle C_\delta:\delta \in S
\rangle,C_\delta \subseteq \delta,\otp(C_\delta) = \mu,S \subseteq \lambda$
\sn
\item[$(iii)$]  $J =: \{A \subseteq \lambda$: for some
club $E$ of $\lambda,\delta \in S \cap A \Rightarrow C_\delta
\nsubseteq E\}$
\sn
\item[$(iv)$]  $\lambda \notin J$ and $\alpha < \lambda \Rightarrow
\lambda > |\{C_\delta \cap \alpha:\alpha \in \nacc(C_\delta),\delta \in S\}|$,
\end{enumerate}
\sn
\item[(d)]   $\mathbf U_J(\lambda_1) < \lambda_2$
\sn
\item[(e)]   for some ${\cP}_1,{\cP}_2$ we have
\sn
\begin{enumerate}
\item[$(i)$]  $\cP_1 \subseteq [\lambda_1]^\kappa,
{\cP}_2 \subseteq [\sigma]^\kappa$
\sn
\item[$(ii)$]   if $g:\sigma \rightarrow \lambda_1$ is one to
one \then \, for some $X \in {\cP}_2$, we have
$\{g(i):i \in X\} \in {\cP}_1$
\sn
\item[$(iii)$]  $|{\cP}_1| < \lambda_2$
\sn
\item[$(iv)$]  $|{\cP}_2| \le \lambda_1$.
\end{enumerate}
\end{enumerate}
\end{claim}

\begin{remark}
\label{c4}
1) We can in \ref{c3} replace clause (a) by
\mn
\begin{enumerate}
\item[$(a)'$]   ${\gk}$ is an a.e.c. which has the $\varphi$-oak
property, see Definition \ref{b3} and $\LST(\gk) \le \lambda_2$.
\end{enumerate}
\mn
2) The proof also gives $\univ(\lambda,\lambda_1,T) \ge \lambda_2$.
\end{remark}

\noindent
Recall
\begin{definition}
\label{c5}
Assume $T,\lambda,\mu,S,\bar C$ are as in Claim \ref{c3}, see (a),(c).

\noindent
1) For $\bar N = \langle N_\gamma:\gamma < 
\lambda \rangle$ an elementary-increasing continuous
sequence of models of $T$ of size $< \lambda$ and for $a,c \in
N_\lambda = \bigcup\limits_{\gamma < \lambda} N_\gamma$ and $\delta \in S$, we
let $\inv_{\varphi,\bar N}(c,C_\delta,a) = \{\zeta < \mu$: there is $b
\in N_{\alpha(\delta,\zeta +2)} \setminus
N_{\alpha(\delta,\zeta +1)}$ such that $N_\lambda \models \varphi[a,b,c])\}$. 

\noindent
2) For $\delta,\bar N$ as above and a set $A \subseteq N_\lambda$, let
 $\inv^A_{\varphi,\bar N}(c,C_\delta) = \bigcup\{\inv_{\varphi,\bar N}
(c,C_\delta,a):a \in A\}$.
\end{definition}

\begin{PROOF}{\ref{c3}}
\medskip

\noindent
\underline{Step A}:  Assume toward a contradiction $\theta =: 
\univ(\lambda_1,T) < \lambda_2$, so let $\langle N^*_j:j < \theta
\rangle$ exemplify this and let $\theta_1 = \theta + |\cP_1| + |\cP_2| +
|T| + {\mathbf U}_J(\lambda_1)$ hence $\theta_1 < \lambda_2$.

Without loss of generality the universe of $N^*_j$ is $\lambda_1$.
\medskip

\noindent
\underline{Step B}:  By the definition of $\mathbf U_J(\lambda_1)$ there is 
${\cA}$ such that:
\mn
\begin{enumerate}
\item[$(a)$]   ${\cA} \subseteq [\lambda_1]^\lambda$
\sn
\item[$(b)$]  $|{\cA}| \le \mathbf U_J(\lambda_1)$
\sn
\item[$(c)$]   if $f:\lambda \rightarrow \lambda_1$ then for some $A
\in {\cA}$ we have $\{\delta \in S:f(\delta) \in A\} \ne \emptyset$
mod $J$.
\end{enumerate}
\mn
For each $X \in {\cP}_1,j < \theta$ and $A \in {\cA}$ let
$M_{j,X,A}$ be an elementary submodel of $N^*_j$ of cardinality
$\lambda$ which includes $X \cup A \subseteq \lambda_1$, and let $\bar
M_{j,X,A} = \langle M_{j,X,A,\varepsilon}:\varepsilon < \lambda
\rangle$ be a filtration of $M_{j,X,A}$.

Lastly, consider

\[
{\cB} = \{\text{inv}^X_{{\bar M}_{j,X,A}}(a,C_\delta):j < \theta,X
\in {\cP}_1,A \in {\cA},\delta \in S \text{ and } a \in
M_{j,X,A}\}.
\]

\noindent
\underline{Step C}:  Easily we have $|{\cB}| \le \theta_1 < \lambda_2$
and $\cB \subseteq [\mu]^\kappa$,
hence there is $B^* \in [\mu]^\kappa \backslash {\cB}$.  \Wilog \,
$\otp(B) = \kappa$, each $\alpha \in B$ is a successor ordinal.

[Why?  Let $h:\mu \rightarrow \mu$ be such that $(\forall \beta <
\mu)(\exists^\mu \alpha < \mu)(h(\alpha +1) = \beta)$ and let $\cB' =
\{\{h(\beta):\beta \in B\}:B \in \cB\}$, so $|\cB'| \le |\cB|$ hence
we can choose $B' \in [\mu]^\kappa \backslash \cB'$.  Let $\langle
\beta_i:i < \kappa \rangle$ list $B'$ and by induction on $i < \kappa$
choose $\alpha_i < \mu$ which is $> \bigcup\limits_{j<i} \alpha_j$ and
satisfies $h(\alpha_i +1) = \beta_i$.
So $\{\alpha_i+1:i < \kappa\}$ is as required.]

Let $\langle \alpha^*_i:i < \kappa\rangle$ list $B$ in increasing
order.  For $\delta \in S$ and $i < \kappa$ let $\alpha_{\delta,i}$ be the
$\alpha^*_i$-th member of $C_\delta$.  Now for $\delta \in S$ and $j
\le \kappa$ let $\nu_\delta = \langle \alpha_{\delta,i}:i < \kappa\rangle$.

Now let
$M^*$ be a $\lambda^+$-saturated model of $T$, in which
$a_i(i < \sigma),b_\eta(\eta \in {}^{\kappa >}(\lambda_2)),c_\nu \, (\nu \in
{}^\kappa(\lambda_2)),\varphi$ are as in the definition of the oak
property and for each $Y \in {\cP}_2$, choose $\langle
N_{Y,\varepsilon}:\varepsilon < \lambda \rangle,\langle 
\bar c_{Y,\varepsilon,\delta}:\delta \in S \rangle$ such that:
\mn
\begin{enumerate}
\item[$(a)$]  $N_{Y,\varepsilon}$ is increasing continuous with
  $\varepsilon$
\sn
\item[$(b)$]  $N_{Y,\varepsilon}$ has cardinality $< \lambda$ for 
$\varepsilon < \lambda$
\sn
\item[$(c)$]  $\bar a_j \in N_{Y,0}$ for $j \in Y$
\sn
\item[$(d)$]  $b_{\nu_\delta \rest (i+1)} \in N_{Y,\nu_\delta(i)+1}$
  for $\delta \in S,i < \kappa$
\sn
\item[$(e)$]  $\bar c_{\nu_\delta} \in N_{Y,\delta +1}$ for $\delta
  \in S$
\sn
\item[$(f)$]  if $i < \kappa,j \in Y,\otp(j \cap Y) = i$ and $\delta
  \in S$ then $N_{Y,\delta +1} \models \varphi[\bar a_j,\bar
  b_{\nu_\delta \rest (i+1)},\bar c_{\nu_\delta}]$.
\end{enumerate}
\mn
As $|\cP_2| \le \lambda_1$ we can choose 
$N \prec M^*,\|N\| = \lambda_1$ such that $\{a_i:i < \sigma\} \cup
\cup \{N_{Y,\varepsilon}:Y \in {\cP}_2,\varepsilon < \lambda\}
\subseteq N$.
\medskip

\noindent
\underline{Step D}:  By our choice of $\langle N^*_j:j < \theta \rangle$,
there is $j(*) < \theta$ and elementary embedding $f:N \rightarrow
N^*_j$.  By an assumption there is $Y \in {\cP}_2$ such that
$X := \{f(\bar a_i):i \in Y\} \in {\cP}_1$.  Also by the choice of ${\cA}$ 
there is $A \in {\cA}$ such that $\{\delta \in S:f(\bar c_{Y,\delta}) 
\in A\} \ne \emptyset \mod J$.

Now we can finish (note that we use here again the last clause in the
definition of the oak property).
\end{PROOF}

\begin{definition}
\label{c7}
1) The formula $\varphi(\bar x,\bar y,\bar z)$
has the weak oak property in $T$ (a first order complete theory) \when
\,: as in Definition \ref{b3} omitting clause (c), (equivalently 
in \cite[1.8]{Sh:710}). 

\noindent
2) A complete first order theory $T$ has the weak oak property \when
\, some $\varphi(\bar x,\bar y,\bar z)$ has it in $T$.

\noindent
3) For non-complete first order property $T$ (or class $\gk =
(K_{\gk},\le_{\gk}))$ we mean $\varphi$ is quantifier-free.
\end{definition}

\noindent
The weak oak property is sufficient for many results on
$\univ(\lambda,T) \ge \lambda_2$ because of
\begin{claim}
\label{c9}
Assume
\mn
\begin{enumerate}
\item[$(a)$]  $T$ has the weak oak property, $|T| \le \lambda = \cf(\lambda)$
\sn
\item[$(b)$]   $\bar C = \langle C_\delta:\delta \in S \rangle,J$
are as in clause (c) of \ref{c3}
\sn
\item[$(c)$]   $\kappa = \cf(\mu) < \sigma < \mu < \lambda =
  \cf(\lambda)$ and $\cP \subseteq \{u \subseteq \sigma:\otp(u) =
  \kappa\}$ has cardinality $\le \lambda$.
\end{enumerate}
\mn
\Then \, for each $B^* \subseteq \mu$ of order type $\kappa,T$ has 
a model $N^*$ of cardinality $\lambda$ and sequence 
$\langle a_i:i < \sigma \rangle$ of members of $N^*$ satisfying:
\mn
\begin{enumerate}
\item[$\circledast$]   if $N$ is a model of $T$ of cardinality
$\lambda$ with filtration $\bar N = \langle N_\alpha:\alpha < \lambda
\rangle$ and $f$ is an elementary embedding of $N^*$ into $N$ \then \,
for every increasing sequence $\bar\varepsilon = \langle
\varepsilon(i):i < \kappa\rangle$ enumerating in increasing order some
$u \in \cP$ we have

\begin{equation*}
\begin{array}{clcr}
\{\delta \in S:&\text{ for some } a \in N^* \text{ we have} \\
  &B^* = \inv^{\{f(a_{\varepsilon(i)}:i < \kappa\}}_{\varphi,\bar N}
(C_\delta,a) = S \mod J.
\end{array}
\end{equation*}
\end{enumerate}
\end{claim}

\begin{PROOF}{\ref{c9}}
\Wilog \, some $\varphi = \varphi(x,y,z)$ witness $T$ has the weak oak
property (as we can replace $T$ by such $T'$ with $\univ(\lambda,T) =
\univ(\lambda,T')$.  

As usual, there is $N^* \models T$ with filtration $\bar
N^* = \langle N^*_i:i < \lambda \rangle$ and 
$I \subseteq {}^{\kappa >} \lambda$ of cardinality $\lambda,\langle
a_i:i < \kappa\rangle,\langle b_\eta:\eta \in {\cT} \rangle$ and 
$\nu_\delta \in
{}^\kappa(C_\delta) \cap \lim_\kappa(T)$ for $\delta \in S$ and
$\langle c_{\nu_\delta}:\delta \in S \rangle$ such that
\mn
\begin{enumerate}
\item[$(a)$]  $\langle a_i:i < \kappa \rangle,\langle b_\eta:\eta
\in {\cT} \rangle,\langle c_{\nu_\delta}:\delta \in S \rangle$ are
as in the Definition \ref{c7}
\sn
\item[$(b)$]   $\otp(\nu_\delta(i) \cap C_\delta) =$ (the $i$-th member
of $B^*) +1$.
\end{enumerate}
\mn
So let $N,\langle N_\varepsilon:\varepsilon < \lambda \rangle,f$ be as
in the assumption of $\circledast$ of the claim.  
Without loss of generality the universes of $N^*$ and of $N$ are $\lambda$.

Let

\[
E_* = \{\delta < \lambda:\delta \text{ limit}, f''(\delta) =
\delta,|N_\delta| = \delta = |N^*_\delta| \text{ and }
(N_\delta,N^*_\delta,f) \prec (N,N^*,f)\}
\]

\mn
it is a club of $\lambda$.  For each $i < \sigma$ let

\begin{equation*}
\begin{array}{clcr}
W_i = \{\alpha:&\text{for some } \delta \in S,\alpha \in C_\delta
\subseteq E,\nu_\delta(i) > \alpha, \\
  &\text{ but } \varphi(f(a_i),y,f(c_{\nu_\delta})) \text{ is
satisfied (in } N) \\
  &\text{ by some } b \in N_\alpha\}.
\end{array}
\end{equation*}

\mn
Now
\mn
\begin{enumerate}
\item[$\circledast$]  $W_i$ is not stationary.
\end{enumerate}
\mn
[Why?  Otherwise let ${\gB} \prec ({\cH}(\lambda^+),\in,<^*)$ be such
that $\bar N,\bar N^*,a_i$ (and even $\langle a_j:j < \sigma\rangle$
and $\cP$ but not used) and $\langle b_\eta:\eta \in {\cT} \rangle,\langle
c_{\nu_\delta}:\delta \in S \rangle$ belong to $\gB$ and ${\gB} \cap
\lambda = \alpha \in W_i$ and assume $b \in {\gB} \cap \alpha,N
\models \varphi(f(a_i),b,f(c_{\nu_\delta}))$.  So there is $\delta(*)
\in S \cap \delta$ such that $N \models
\varphi[f(a_1),b,f(c_{\nu_{\delta(*)}})$.  But $\nu_\delta(i) \ge \alpha >
\nu_{\delta(*)}(i)$ hence
$\varphi(a_i,y,c_{\nu_\delta}),\varphi(a_i,y,c_{\nu_{\delta'}})$ are
incompatible (in $N^*$) hence their images by $f$ are incompatible in
$N$ by $b$ satisfies both, a contradiction, so $W_i$ is not stationary.]

So there is a club $E^*$ of $\lambda$ included in $E_*$ and
disjoint to $W_i$ for each $i < \sigma$.  So there is $\delta \in S$
such that $C_\delta \subseteq E^*$ and we get contradiction as earlier.
\end{PROOF}

\begin{question}
\label{c11}
Can we combine \ref{c3}, \ref{c9}? 

(For many singular $\lambda_1$'s, certainly yes).
\end{question}
\newpage


\bibliographystyle{amsalpha}
\bibliography{shlhetal}

\end{document}